\documentclass[11pt]{amsart}

\usepackage{amssymb}
\usepackage{amsmath}
\usepackage{amsthm, geometry, graphicx}
\usepackage{array}
\usepackage{tikz}

\newtheorem{thm}{Theorem}
\newtheorem*{thm*}{Theorem}

\author{K. Blum}
\title{Bounds on the Number of Graphical Partitions}
\begin{document}

\maketitle
\begin{abstract}
    We narrow in on the number of graphical partitions for which there is no known generating function by manipulating the well known generating function for Frobenius partitions. 
\end{abstract}
\section{Introduction}
A graphical partition of $n$ is a partition of $n$ where each element corresponds to the degree of a vertex in a simple undirected graph. For example, $\{4,2,2,2,2\}$ is a graphical partition of 12 because it is the degree sequence for a realizable simple undirected graph with 6 edges. The set of graphical partitions is a subset of the set of all integer partitions. That is, all graphical partitions are integer partitions, but not all integer partitions are graphical.

Let $g(n)$ be the number of graphical partitions of $n$ and let $p(n)$ be the number of integer partitions of $n$. Currently, there is no known generating function for $g(n)$ and there are only a handful of ways to determine $g(n)$ for any given $n$. One of the most accurate but least efficient ways is to test each individual integer partition of $n$ against the Erd\"os-Gallai theorem \cite{ErdosandGallai} and eliminate the sequences which are not graphical. 

Consider the partition $\lambda=\{d_1, d_2, d_3, \cdots, d_{\ell}\}$. By looking at the Ferrer's diagram of $\lambda$ we can construct the corresponding Frobenius partition
\begin{equation*}{}
\Big( \begin{array}{cccc}
     a_1& a_2 & \cdots & a_k  \\
     b_1& b_2 & \cdots & b_k
\end{array}{}\Big)
\end{equation*}
wherein \[a_m=d_m-m\]
and
\[b_m=\frac{m(m-1)}{2}+\sum_{i=m+1}^{\ell}\min\{m,d_i\}-\sum_{i=1}^{m-1}b_i.\]

We define $f(n)$ to be the number of Frobenius partitions of $n$ wherein the sum of the entries on the top row plus the number of columns is less than or equal to the sum of the entries on the bottom row and modify the generating function for Frobenius partitions to enumerate $f(n)$. That is, we define $f(n)$ to be the number of partitions of weight $n$ wherein \[k+\sum_{i=1}^{k}a_i\leq \sum_{i=1}^{k}b_i.\] 
Next, let the $m$-th successive rank be defined as $a_m-b_m$. We define $f'(n)$ to be the number of Frobenius partitions that satisfy the condition that each successive rank is less than or equal to $-1$. That is, we define $f'(n)$ to be the number of Frobenius partitions of weight $n$ wherein \[a_m-b_m \leq -1\] for each $1\leq m \leq k$.  

We show that $f'(n)\leq g(n)\leq f(n)\leq p(n)$. % and investigate limits analogous to those established by Erd\"os and Richmond in \cite{ErdosandRichmond}. 

\section{Generating Functions for $f(n)$ and $f'(n)$}

%In combinatorics, a generating function is a function which generates a formal power series who has, as its coefficients, a particular sequence. Formally stated, it is a formal power series 
%\begin{equation*}
%    f(q)=\sum_{n=0}^{\infty}a_nq^n
%\end{equation*}{}
%whose coefficients are the sequence $\{a_0,a_1,a_2,...\}$.

In \cite{Andrews} Andrews establishes a generating function that enumerates both $p(n)$ and the total number of Frobenius partitions of weight $n$. He shows that the coefficient of $z^0$ in
\begin{equation*}\label{generating function}
     \prod_{i=0}^{\infty}(1+zq^{i+1})(1+z^{-1}q^i)
\end{equation*}{} is an infinite series in $q$ wherein the coefficient of $q^n$ is $p(n)$, or alternatively, wherein the coefficient of $q^n$ is the number of Frobenius partitions of weight $n$. 

We can modify the established generating function for the number of Frobenius partitions to enumerate $f(n)$ by separating the sum of the entries on the bottom row from the sum of entries on the top row plus the number of columns. Then, we have the generating function 
\begin{equation*}
     \prod_{i=0}^{\infty}(1+zq^{i+1})(1+z^{-1}x^i)
\end{equation*}{}
where the coefficient of $z^0$ is a series in $q$ and $x$, and the coefficient of $q^mx^r$ where $m+r=n$ and $m\leq r$ is $f(n)$. 

In \cite{Andrews} and \cite{Bressoud} Andrews and Bressoud establish the following theorem:  
 \begin{thm*}{}(Andrews-Bressoud)
 Given positive integer $M$ and integral $r$ with $0\leq r \leq M/2$ let $A_{M,r}(n)$ denote the number of partitions into parts not congruent to $0$ or $\pm r$ modulo $M$. Let $B_{M,r}(n)$ denote the number of partitions of $n$ whose successive ranks lie in the interval $[-r+2, M-r-2]$. Then, $A_{M,r}(n)=B_{M,r}(n)$ for all $n$
 \end{thm*}
 
 To enumerate $f'(n)$ we wish to count the number of partitions whose successive ranks lie in the interval $[-n+1,-1]$. Therefore, we can let $M=n+2$ and $r=n+1$ and count the number of partitions where no part is congruent to 0 or $\pm n+1$ modulo $n+2$. Then, since $-(n+1) \equiv 1 \pmod{n+2}$ and no part of a partition of $n$ can be bigger than $n$, the problem of enumerating $f'(n)$ boils down to counting the number of partitions that do not contain 1 as a part. Therefore, the generating function for $f'(n)$ is \[\prod_{i=2}^{\infty}\frac{1}{1-q^i}.\]

%We intend to show that $g(n)$ is a subset of $f(n)$ which is a subset of $p(n)$ and therefore $g(n)\leq f(n) \leq p(n)$. Then, we show that there is a condition on Frobenius partitions counted by $f(n)$ that ensures that the partition is graphical. 

\section{Establishing $f'(n)\leq g(n)\leq f(n)\leq p(n)$}
Since there is a bijection between the set of integer partitions and the set of all Frobenius partitions, we can conclude that the set of Frobenius partitions with the condition that the sum of the entries on the top row plus the number of columns is less than or equal to the sum of the entries on the bottom row is a subset of the set of integer partitions. That is the set of partitions counted by $f(n)$ is a subset of the partitions counted by $p(n)$. Therefore, $f(n)\leq p(n)$. We show that the set of graphical partitions is a subset of the set of Frobenius partitions with the condition that the sum of the entries on the top row plus the number of columns is less than or equal to the sum of the entries on the bottom row. In other words, we show that if a partition is counted by $g(n)$ then it is also counted by $f(n)$ and therefore $g(n)\leq f(n)$.  
\begin{thm}
Let $\lambda=\{d_1,d_2,\dots,d_{\ell}\}$ be a nonincreasing sequence of nonnegative integers such that $d_1+d_2+\cdots+d_{\ell}=2n$ for some $n\in \mathbb{Z^+}$. Furthermore, let
\begin{equation*}{}
\Big( \begin{array}{cccc}
     a_1& a_2 & \cdots & a_k  \\
     b_1& b_2 & \cdots & b_k
\end{array}{}\Big)
\end{equation*}
be the Frobenius partition of $\lambda$.  
If $\lambda$ is a graphical partition of $2n$ then
\begin{equation*}
    k+\sum_{i=1}^{k}a_i\leq \sum_{i=1}^{k}b_i.
\end{equation*}
\end{thm}{}

\begin{proof}
Assume $\lambda$ is a graphical partition. Suppose to the contrary that 
\begin{equation*}
    k+\sum_{i=1}^{k}a_i>\sum_{i=1}^{k}b_i.
\end{equation*}{}
From Erd\"os and Gallai, we must have that
\begin{equation*}
    \sum_{i=1}^{k}d_i\leq k(k-1)+\sum_{i=k+1}^{\ell}\min\{k,d_i\}.
\end{equation*}{}
By looking at the Ferrer's diagram construction of the Frobenius partition, we have that
\begin{equation*}\label{toprow}
   k+\sum_{i=1}^{k}a_i=k-\frac{k(k+1)}{2}+ \sum_{i=1}^{k}d_i
\end{equation*}{}
and
\begin{equation*}
    \sum_{i=1}^{k}b_i=\frac{k(k-1)}{2}+\sum_{i=k+1}^{\ell}\min\{k,d_i\}.
\end{equation*}{}

So, we have that 
\begin{equation*}
    \frac{k(k-1)}{2}+\sum_{i=k+1}^{\ell}\min\{k,d_i\}< k-\frac{k(k+1)}{2}+ \sum_{i=1}^{k}d_i.
\end{equation*}{}
and that 
\begin{equation*}
    k- \frac{k(k+1)}{2}+\sum_{i=1}^{k}d_i\leq k-\frac{k(k+1)}{2}+k(k-1)+\sum_{i=k+1}^{\ell}\min\{k,d_i\}.
\end{equation*}
So, 
\begin{equation*}
    \frac{k(k-1)}{2}<k-\frac{k(k+1)}{2}+k(k-1)
\end{equation*}{}
which implies that
\begin{equation*}
    k^2-k<k^2-k
\end{equation*}{}
and that is a contradiction. Therefore, we conclude that if $\lambda$ is a graphical partition then 
\begin{equation*}
    k+\sum_{i=1}^{k}a_i\leq \sum_{i=1}^{k}b_i.
\end{equation*}{}
\end{proof}{}

 Thus, we have established that the set of partitions counted by $g(n)$ is a subset of of the partitions counted by $f(n)$.
 Next, we'll establish that the set of partitions counted by $f'(n)$ is a a subset of the set of partitions counted by $g(n)$.

\begin{thm}
Let $\lambda=\{d_1,d_2,\dots,d_{\ell}\}$ be a nonincreasing sequence of nonnegative integers such that $d_1+d_2+\cdots+d_{\ell}=2n$ for some $n\in \mathbb{Z^+}$. Furthermore, let
\begin{equation*}{}
\Big( \begin{array}{cccc}
     a_1& a_2 & \cdots & a_k  \\
     b_1& b_2 & \cdots & b_k
\end{array}{}\Big)
\end{equation*}
be the Frobenius partition of $\lambda$. If $a_m-b_m\leq -1$ for all $1\leq m \leq k \leq \ell$ then $\lambda$ is a graphical partition. 
\end{thm}{}

\begin{proof}
Assume 
\[a_m-b_m\leq -1\] for all $1\leq m \leq k \leq \ell$. 

So, 
\[a_m+1\leq b_m.\]

Now, by examining the Ferrer's diagram construction of the Frobenius partition, we have
\[a_m=d_m-m\]
and
\[b_m=\frac{m(m-1)}{2}+\sum_{i=m+1}^{\ell}\min\{m,d_i\}-\sum_{i=1}^{m-1}b_i.\]
So, 
\[d_m-m+1\leq \frac{m(m-1)}{2}+\sum_{i=m+1}^{\ell}\min\{m,d_i\}-\sum_{i=1}^{m-1}b_i\]
We need to show that
\[\sum_{i=1}^{m}d_i\leq m(m-1)+\sum_{i=m+1}^{\ell}\min\{m,d_i\}\]
to complete the proof. 
So, again, we have
\begin{align*}
    d_m-m+1&\leq \frac{m(m-1)}{2}+\sum_{i=m+1}^{\ell}\min\{m,d_i\}-\sum_{i=1}^{m-1}b_i\\
    d_m & \leq m-1+\frac{m(m-1)}{2}+\sum_{i=m+1}^{\ell}\min\{m,d_i\}-\sum_{i=1}^{m-1}b_i
\end{align*}
So, 
\begin{align*}
   \sum_{i=1}^{m}d_i&\leq m-1+\frac{m(m-1)}{2}+\sum_{i=m+1}^{\ell}\min\{m,d_i\}-\sum_{i=1}^{m-1}b_i+\sum_{i=1}^{m-1}d_i\\
    &= m-1+\frac{m(m-1)}{2}+\sum_{i=m+1}^{\ell}\min\{m,d_i\}-\sum_{i=1}^{m-1}b_i+\sum_{i=1}^{m-1}[a_i+i]\\&=m-1+m(m-1)+\sum_{i=m+1}^{\ell}\min\{m,d_i\}-\sum_{i=1}^{m-1}b_i+\sum_{i=1}^{m-1}a_i\\
    &=m-1+m(m-1)+\sum_{i=m+1}^{\ell}\min\{m,d_i\}+\sum_{i=1}^{m-1}[a_i-b_i]\\
\end{align*}{}
Now, since we assume $a_i-b_i\leq -1$ for all $i$, 
\[\sum_{i=1}^{m-1}[a_i-b_i]\leq -(m-1).\] 
So, 
\[\sum_{i=1}^{m}d_i\leq m(m-1)+\sum_{i=m+1}^{\ell}\min\{m,d_i\}\]
Therefore, if $a_m-b_m\leq -1$ for all $1\leq m \leq k \leq \ell$ then
\[\sum_{i=1}^{m}d_i\leq m(m-1)+\sum_{i=m+1}^{\ell}\min\{m,d_i\}\] and $\lambda$ is graphical. 
\end{proof}{}

\section{A Note On Dyson Rank}
The Dyson rank of a partition is defined as the largest part minus the number of parts. We provide a new proof that if a partition is graphical then its Dyson rank is less than or equal to -1. 
\begin{thm}
 Let $\lambda=\{d_1,d_2,..., d_{\ell}\}$ be a sequence of nonincreasing nonnegative integers and such that $d_1+d_2+\cdots+d_{\ell}=2n$. Let \[\Big(\begin{array}{ccccc}
     a_1 & a_2 & \cdots & a_k  \\
    b_1 & b_2 & \cdots & b_k
\end{array}\Big)\]
be the Frobenius partition of $\lambda$. Then the Dyson rank of $\lambda$ is $r=a_1-b_1$ and if $\lambda$ is a graphical partition of $2n$ then $r\leq -1$.
\end{thm}{}

\begin{proof}

Assume $\lambda$ is graphical. 
From Erd\"os and Gallai, 
\[d_1\leq \sum_{i=2}^{\ell}\min\{1,d_i\}\leq \ell-1\]
So, since rank is defined as the largest part of the partition minus the number of parts, and since $d_1$ is the largest part and $\ell$ is the number of parts, then the rank $k$ is $d_1-\ell$ and
\[d_1-\ell\leq -1.\] 
Furthermore, since
\[d_1=a_1+1 \quad{and} \quad \ell=b_1+1\]
we also have that if $\lambda$ is graphical then 
$a_1-b_1\leq -1$. 

\end{proof}{}

\section{Table of Values}
The following is a table of values for the first 20 even values of $n$: 
\begin{center}
\begin{tabular}{|m{1.5cm}|m{1.5cm}|m{1.5cm}|m{1.5cm}|m{1.5cm}|}\hline
$\mathbf{n}$  & $f'(n)$ & $g(n)^{\star}$  & $f(n)$ & $p(n)$ \\ 
\hline 2 & 1 & 1 & 1 & 2 \\ 
\hline 4 & 2 & 2 & 2 & 5 \\
\hline 6 & 4 & 5 & 5 & 11 \\
\hline 8 & 7 & 9 & 9 & 22 \\
\hline 10 & 12 & 17 & 18 & 42 \\
\hline 12 & 21 & 31 & 32 & 77 \\
\hline 14 & 34 & 54 & 57 & 135\\
\hline 16 & 55 & 90 & 95 & 231\\
\hline 18 & 88 & 151 & 162 & 385\\
\hline 20 & 137 & 244 & 264 & 627\\
\hline 22 & 210 & 387 & 418 & 1002\\
\hline 24 & 320 & 607 & 659 & 1575\\
\hline 26 & 478 & 933 & 1016 & 2436\\
\hline 28 & 708 & 1420 & 1555 & 3718 \\
\hline 30 & 1039 & 2136 & 2347 & 5604 \\
\hline 32 & 1507 & 3173 &  3499 & 8349 \\
\hline 34 & 2167 &  4657 & 5152 &  12310\\
\hline 36 & 3094 & 6799 & 7558 & 17977  \\
\hline 38 & 4378 & 9803 & 10914 & 26015  \\
\hline 40 & 6153 & 14048 & 15704 & 37338\\
\hline
\end{tabular}
\end{center}

\section{Acknowledgements}
A big thank you to Dr. Taylor Short, Kenneth Moore and Drew Posh from Grand Valley State University for letting me run with this idea and to Dr. Andrej Dudek from Western Michigan University and Dr. Nicholas Ovenhouse from University of Minnesota for their comments.

\end{document}